\documentclass[12pt, twoside, leqno]{article}

\usepackage{amsmath,amsthm}
\usepackage{amssymb}

\numberwithin{equation}{section}
\newtheorem{thm}{Theorem}[section]
\newtheorem{defn}[thm]{Definition}

\newtheorem{lem}[thm]{Lemma}

\newtheorem{re}{Remark}[section]
\newenvironment{pf}{{\noindent \it \bf Proof:}}{{\hfill$\Box$}}

\begin{document}

\title{A regularity criterion for 3D micropolar fluid flows in terms of one
partial derivative of the velocity}
\author{Sadek Gala \\
{\small Department of Mathematics, University of Mostaganem}\\
{\small Box 227, Mostaganem, Algeria}\\
\&\\
{\small Dipartimento di Mathematica e Informatica}\\
{\small Universit\`{a} di Catania, Viale Andrea Doria, 6, 95125 Catania,
Italy}\\
{\small sadek.gala@gmail.com}\\
and\\
Maria \ Alessandra Ragusa\\
{\small Dipartimento di Mathematica e Informatica}\\
{\small Universit\`{a} di Catania, Viale Andrea Doria, 6, 95125 Catania,
Italy}\\
{\small maragusa@dmi.unict.it}}
\date{}
\maketitle

\begin{abstract}
In this work, we prove a regularity criterion for micropolar fluid flows in
terms of the one partial derivative of the velocity in Morrey-Campanato
space.
\end{abstract}


\noindent {\textit{Mathematics Subject Classification(2000):}\thinspace
\thinspace 35Q35,\thinspace \thinspace 35B65,\thinspace \thinspace 76D05}

\section{Introduction and the main result}

\bigskip In this paper, we consider the following Cauchy problem for the
incompressible micropolar fluid equations in $\mathbb{R}^{3}$ \cite{E}:%
\begin{equation}
\left\{
\begin{array}{l}
\partial _{t}u+\left( u\cdot \nabla \right) u-\Delta u+\nabla \pi -\nabla
\times \omega =0, \\
\partial _{t}\omega -\Delta \omega -\nabla div\omega +2\omega +u\cdot \nabla
\omega -\nabla \times u=0, \\
\nabla \cdot u=0, \\
u(x,0)=u_{0}(x),\text{ }\omega (x,0)=\omega _{0}(x),%
\end{array}%
\right.  \label{eq1}
\end{equation}%
where $u$, $\omega $ and $\pi $ denote the unknown velocity vector field,
the micro-rotational velocity and the unknown scalar pressure of the fluid
at the point $\left( x,t\right) \in \mathbb{R}^{3}\times (0,T)$,
respectively, while $u_{0},\omega _{0}$ are given initial data with $\nabla
\cdot u_{0}=0$\ in the sense of distributions.

When the micro-rotation effects are neglected or $\omega =0$, the micropolar
fluid flows (\ref{eq1}) reduce to the incompressible Navier-Stokes flows
(see, for example, \cite{La,Te}). Much effort has been devoted to establish
the global existence and uniqueness of smooth solutions to the
Navier--Stokes equations. Different criteria for regularity of the weak
solutions have been proposed. The Prodi--Serrin conditions (see \cite{Pro,
Ser}) shows that any solution $u$ for the 3D Navier-Stokes equations
satisfying
\begin{equation}
u\in L^{p}(0,T;L^{q}(\mathbb{R}^{3}))\text{ \ with \ }\frac{2}{p}+\frac{3}{q}%
\leq 1\text{ \ and \ }3\leq q\leq \infty ,  \label{eq1.5}
\end{equation}%
is regular. Notice that the limiting case $u\in L^{\infty }(0,T;L^{3}(%
\mathbb{R}^{3}))$ was covered by Escauriaza et al. \cite{Esc} in 2003. Later
on, Beir\~{a}o da Veiga \cite{Bei} established another regularity criterion
by replacing (\ref{eq1.5}) with the following condition:
\begin{equation}
\nabla u\in L^{\beta }\left( (0,T;L^{\alpha }(\mathbb{R}^{3})\right) \text{
\ with \ }\frac{3}{\alpha }+\frac{2}{\beta }\leq 2\text{ \ and \ }\frac{3}{2}%
<\alpha \leq \infty .  \label{eq1.6}
\end{equation}%
In 2004, Penel and Pokorn\'{y} \cite{PP} obtained a different type
regularity criterion, which says that if
\begin{equation}
\partial _{3}u\in L^{\beta }\left( (0,T;L^{\alpha }(\mathbb{R}^{3})\right)
\text{ \ \ with \ }\frac{3}{\alpha }+\frac{2}{\beta }\leq 1\text{ \ and }%
2\leq \alpha \leq \infty ,  \label{eq1.7}
\end{equation}%
then the solution $u$ to the Navier-Stokes equations is regular. The same
result can be found in \cite{Z2}. Penel and Pokorn\'{y}'s work has been
improved by some other authors, (see e.g., \cite{Cao, KZ} and the references
cited therein). It was already known that if one component of the velocity
is bounded in a suitable space, then the solution is smooth (see Penel and
Pokorn\'{y} \cite{PP}, Zhou \cite{Z1, Z2, Z3, Z4}). Some of these regularity
criteria can be extended to the 3D MHD equations by making the assumptions
on both $u$ and $b$ (\cite{CKS}). Moreover, He and Xin in \cite{HX} derived
some regularity criteria for the 3D MHD equations only in terms of the
velocity field $u$, and they proved that if u satisfies either (\ref{eq1.5})
or (\ref{eq1.6}), then the solution is regular. Recently, Cao and Wu \cite%
{CW} proved that the condition
\begin{equation}
\partial _{3}u\in L^{\beta }\left( (0,T;L^{\alpha }(\mathbb{R}^{3})\right)
\text{ \ \ with \ }\frac{3}{\alpha }+\frac{2}{\beta }\leq \frac{3}{2}\text{
\ and }\alpha >3,  \label{eq1.8}
\end{equation}%
also implies regularity of the solution $(u,b)$ to the 3D MHD equations.
Later, Jia and Zhou \cite{XZ1, XZ2, XZ4} showed that if
\begin{equation}
\partial _{3}u\in L^{\beta }\left( (0,T;L^{\alpha }(\mathbb{R}^{3})\right)
\text{ \ \ with \ }\frac{3}{\alpha }+\frac{2}{\beta }=\frac{3}{4}+\frac{1}{%
\alpha }\text{ \ and }\alpha >2,  \label{eq1.9}
\end{equation}%
then the solution is regular. For more interesting component reduction
results of the regularity criterion, we refer to e.g. \cite{XZ3, Z1, Z2, Z3,
Z4}.

Inspired by the above-mentioned works on regularity criteria of
Navier-Stokes and MHD equations, particularly those of Penel and Pokorn\'{y}
\cite{PP}, Cao and Wu \cite{CW} and Jia and Zhou \cite{XZ1, XZ2, XZ3, XZ4},
we want to investigate a similar problem for the micropolar fluid flows (\ref%
{eq1}). Very recently, Jia et al. \cite{Jia} proved the following regularity
criterion%
\begin{equation*}
\partial _{3}u\in L^{\beta }\left( (0,T;L^{\alpha ,\infty }(\mathbb{R}%
^{3})\right) \text{ \ \ with \ }\frac{3}{\alpha }+\frac{2}{\beta }=1\text{ \
and }3<\alpha \leq \infty .
\end{equation*}%
Here $L^{\alpha ,\infty }$ is the Lorentz space.

The purpose of this work is to improve the result in \cite{Jia} and to prove
that if the derivative of the velocity in one direction belongs to $L^{\frac{%
2}{1-r}}\left( 0,T,\overset{.}{\mathcal{M}}_{2,\frac{3}{r}}(\mathbb{R}%
^{3})\right) $ \ with \ $0<r<1$, then the weak solution actually is regular
and unique. This work is motivated by the recent results \cite{XZ1}-\cite{Z4}
on the Navier-Stokes equations and MHD equations.

\section{Preliminaries and main result}

Now, we recall the definition and some properties of the space that will be
useful in the sequel. These spaces play an important role in studying the
regularity of solutions to partial differential equations; see e.g. \cite{GL}
and references therein.

\begin{defn}
For $0\leq r<\frac{3}{2}$, the space $\overset{.}{X}_{r}$ is defined as the
space of $f(x)\in L_{loc}^{2}\left( \mathbb{R}^{3}\right) $ such that
\begin{equation*}
\left\Vert f\right\Vert _{\overset{.}{X}_{r}}=\underset{\left\Vert
g\right\Vert _{\overset{.}{H}^{r}}\leq 1}{\sup }\left\Vert fg\right\Vert
_{L^{2}}<\infty .
\end{equation*}%
where we denote by $\overset{.}{H}^{r}\left( \mathbb{R}^{3}\right) $ the
completion of the space $C_{0}^{\infty }\left( \mathbb{R}^{3}\right) $ with
respect to the norm $\left\Vert u\right\Vert _{\overset{.}{H}%
^{r}}=\left\Vert \left( -\Delta \right) ^{\frac{r}{2}}u\right\Vert _{L^{2}}.$
\end{defn}

We have the homogeneity properties : $\forall x_{0}\in \mathbb{R}^{3}$%
\begin{align*}
\left\Vert f(.+x_{0})\right\Vert _{\overset{.}{X}_{r}}& =\left\Vert
f\right\Vert _{\overset{.}{X}_{r}} \\
\left\Vert f(\lambda .)\right\Vert _{\overset{.}{X}_{r}}& =\frac{1}{\lambda
^{r}}\left\Vert f\right\Vert _{\overset{.}{X}_{r}}\text{, \ \ }\lambda >0.
\end{align*}%
The following imbedding :%
\begin{equation*}
L^{\frac{3}{r}}\subset \overset{.}{X}_{r},\text{ \ \ \ \ }0\leq r<\frac{3}{2}
\end{equation*}%
holds.

Now we recall the definition of Morrey-Campanato spaces (see e.g. \cite{Kat}%
):

\begin{defn}
For $1<p\leq q\leq +\infty $, the Morrey-Campanato space $\overset{\cdot }{%
\mathcal{M}}_{p,q}$ is defined by :%
\begin{equation}
\overset{\cdot }{\mathcal{M}}_{p,q}=\left\{ f\in L_{loc}^{p}\left( \mathbb{R}%
^{3}\right) :\text{ }\left\Vert f\right\Vert _{\overset{\cdot }{\mathcal{M}}%
_{p,q}}=\underset{x\in \mathbb{R}^{3\text{ \ }}}{\sup }\underset{R>0}{\sup }%
R^{3/q-3/p}\left\Vert f\right\Vert _{L^{p}(B(x,R))}<\infty \right\} .
\label{eq1.13}
\end{equation}
\end{defn}

It is easy to check the following :%
\begin{equation*}
\left\Vert f(\lambda .)\right\Vert _{\overset{\cdot }{\mathcal{M}}_{p,q}}=%
\frac{1}{\lambda ^{\frac{3}{q}}}\left\Vert f\right\Vert _{\overset{\cdot }{%
\mathcal{M}}_{p,q}}\text{, \ \ \ \ \ }\lambda >0.
\end{equation*}

We have the following comparison between Lorentz spaces and Morrey-Campanato
spaces : for $p\geq 2$%
\begin{equation*}
L^{\frac{3}{r}}\left( \mathbb{R}^{3}\right) \subset \text{ }L^{\frac{3}{r}%
,\infty }\left( \mathbb{R}^{3}\right) \subset \overset{.}{\mathcal{M}}_{p,%
\frac{3}{r}}\left( \mathbb{R}^{3}\right) .
\end{equation*}%
Other useful comparison are contained in \cite{JOTA}, \cite{HERZ} and \cite{AMC}. The relation
\begin{equation*}
L^{\frac{3}{r},\infty }\left( \mathbb{R}^{3}\right) \subset \overset{.}{%
\mathcal{M}}_{p,\frac{3}{r}}\left( \mathbb{R}^{3}\right)
\end{equation*}%
is shown as follows.\ Let $f\in L^{\frac{3}{r},\infty }\left( \mathbb{R}%
^{3}\right) $.\ Then
\begin{eqnarray*}
\left\Vert f\right\Vert _{\overset{.}{\mathcal{M}}_{p,\frac{3}{r}}} &\leq &%
\underset{E}{\sup }\left\vert E\right\vert ^{\frac{r}{3}-\frac{1}{2}}\left(
\int\limits_{E}\left\vert f(y)\right\vert ^{p}dy\right) ^{\frac{1}{p}} \\
&=&\left( \underset{E}{\sup }\left\vert E\right\vert ^{\frac{pr}{3}%
-1}\int\limits_{E}\left\vert f(y)\right\vert ^{p}dy\right) ^{\frac{1}{p}}%
\text{ } \\
&\cong &\left( \underset{R>0}{\sup }R\left\vert \left\{ x\in \mathbb{R}%
^{3}:\left\vert f(y)\right\vert ^{p}>R\right\} \right\vert ^{\frac{pr}{3}%
}\right) ^{\frac{1}{p}} \\
&=&\underset{R>0}{\sup }R\left\vert \left\{ x\in \mathbb{R}^{p}:\left\vert
f(y)\right\vert >R\right\} \right\vert ^{\frac{r}{3}} \\
&\cong &\left\Vert f\right\Vert _{L^{\frac{3}{r},\infty }}.
\end{eqnarray*}

For $0<r<1$, we use the fact that
\begin{equation*}
L^{2}\cap \overset{.}{H}^{1}\subset \overset{.}{B}_{2,1}^{r}\subset \overset{%
.}{H}^{r}.
\end{equation*}%
Thus we can replace\ the space $\overset{.}{X}_{r}$ by the pointwise
multipliers from Besov space $\overset{.}{B}_{2,1}^{r}$ to $L^{2}$.\ Then we
have the following lemma given in \cite{Lem1}.

\begin{lem}
\label{lem 2}For $0\leq r<\frac{3}{2}$, the space $\overset{.}{Z}_{r}$ is
defined as the space of $f(x)\in L_{loc}^{2}\left( \mathbb{R}^{3}\right) $
such that
\begin{equation*}
\left\Vert f\right\Vert _{\overset{.}{Z}_{r}}=\underset{\left\Vert
g\right\Vert _{\overset{.}{B}_{2,1}^{r}}\leq 1}{\sup }\left\Vert
fg\right\Vert _{L^{2}}<\infty .
\end{equation*}%
Then $f\in \overset{\cdot }{\mathcal{M}}_{2,\frac{3}{r}}$ if and only if $%
f\in \overset{.}{Z}_{r}$ with equivalence of norms.
\end{lem}

To prove our main result, we need the following lemma due to \cite{MO} (see
also \cite{ZG}).

\begin{lem}
\label{lem3}For $0<r<1$, we have
\begin{equation*}
\left\Vert f\right\Vert _{\overset{.}{B}_{2,1}^{r}}\leq C\left\Vert
f\right\Vert _{L^{2}}^{1-r}\left\Vert \nabla f\right\Vert _{L^{2}}^{r}.
\end{equation*}
\end{lem}

Additionally, for $2<p\leq \frac{3}{r}$ and $0\leq r<\frac{3}{2}$, we have
the following inclusion relations (\cite{Lem}, \cite{Lem1}) :
\begin{equation*}
\overset{.}{\mathcal{M}}_{p,\frac{3}{r}}\left( \mathbb{R}^{3}\right) \subset
\overset{.}{X}_{r}\left( \mathbb{R}^{3}\right) \subset \overset{.}{\mathcal{M%
}}_{2,\frac{3}{r}}\left( \mathbb{R}^{3}\right) =\overset{.}{Z}_{r}\left(
\mathbb{R}^{3}\right) \text{.}
\end{equation*}%
The relation
\begin{equation*}
\overset{.}{X}_{r}\left( \mathbb{R}^{3}\right) \subset \overset{.}{\mathcal{M%
}}_{2,\frac{3}{r}}\left( \mathbb{R}^{3}\right)
\end{equation*}%
is shown as follows. Let $f\in \overset{.}{X}_{r}\left( \mathbb{R}%
^{3}\right) $, $0<R\leq 1$ , $x_{0}\in \mathbb{R}^{3}$ and $\phi \in
C_{0}^{\infty }\left( \mathbb{R}^{3}\right) $, $\phi \equiv 1$ on $B(\frac{%
x_{0}}{R},1)$. We have
\begin{align*}
R^{r-\frac{3}{2}}\left( \int_{\left\vert x-x_{0}\right\vert \leq
R}\left\vert f(x)\right\vert ^{2}dx\right) ^{1/2}& =R^{r}\left(
\int_{\left\vert y-\frac{x_{0}}{R}\right\vert \leq 1}\left\vert
f(Ry)\right\vert ^{2}dy\right) ^{1/2} \\
& \leq R^{r}\left( \int_{y\in \mathbb{R}^{3}}\left\vert f(Ry)\phi
(y)\right\vert ^{2}dy\right) ^{1/2} \\
& \leq R^{r}\left\Vert f(R.)\right\Vert _{\overset{.}{X}_{r}}\left\Vert \phi
\right\Vert _{H^{r}} \\
& \leq \left\Vert f\right\Vert _{\overset{.}{X}_{r}}\left\Vert \phi
\right\Vert _{H^{r}} \\
& \leq C\left\Vert f\right\Vert _{\overset{.}{X}_{r}}.
\end{align*}

Before stating our result, let us recall the definition of Leray--Hopf weak
solution.

\begin{defn}[\protect\cite{Lu}]
\label{def1}Let $\left( u_{0},\omega _{0}\right) \in L^{2}\left( \mathbb{R}%
^{3}\right) $ and $\nabla \cdot u_{0}=0$. A measurable function $\left(
u(x,t),\omega (x,t)\right) $ is called a weak solution to the $3D$
micropolar flows equations (\ref{eq1}) on $\left( 0,T\right) $ if $\left(
u,\omega \right) $ satisfies the following properties

\begin{description}
\item[(1)] $u,\omega \in L^{\infty }\left( \left( 0,T\right) ;L^{2}\left(
\mathbb{R}^{3}\right) \right) \cap L^{2}\left( \left( 0,T\right)
;H^{1}\left( \mathbb{R}^{3}\right) \right) $ for all $T>0$;

\item[(2)] $\left( u(x,t),\omega (x,t)\right) $ verifies (\ref{eq1}) in the
sense of distribution;

\item[(3)] The energy inequality
\begin{flalign*}
&\| u\|_{L^{2}}^{2}+\| w\| _{L^{2}}^{2}+2\int\limits_{0}^{t}\left(
\| \nabla u\| _{L^{2}}^{2}+\| \nabla w\| _{L^{2}}^{2}\right)
ds+2\int\limits_{0}^{t}\| \nabla \mathbf{\cdot}w\|
_{L^{2}}^{2}ds +2\int\limits_{0}^{t}\| w\| _{L^{2}}^{2}ds \\
& \leq \| u_{0}\| _{L^{2}}^{2}+\| w_{0}\| _{L^{2}}^{2},
~\mbox{for}~~ 0<t\leq T.&
\end{flalign*}
\end{description}
\end{defn}

By a strong solution we mean a weak solution $\left( u,\omega \right) $ such
that
\begin{equation*}
u,\omega \in L^{\infty }\left( \left( 0,T\right) ;H^{1}\left( \mathbb{R}%
^{3}\right) \right) \cap L^{2}\left( \left( 0,T\right) ;H^{2}\left( \mathbb{R%
}^{3}\right) \right) .
\end{equation*}%
It is well known that strong solutions are regular (say, classical) and
unique in the class of weak solutions.

More precisely, we will prove

\begin{thm}
\label{th1}Suppose that $(u_{0},\omega _{0})\in H^{1}(\mathbb{R}^{3})$ with $%
\nabla \mathbf{.}u_{0}=0$ in $\mathbb{R}^{3}$. If the velocity $u$ satisfies
\begin{equation}
\partial _{3}u\in L^{\frac{2}{1-r}}\left( 0,T,\overset{.}{\mathcal{M}}_{2,%
\frac{3}{r}}(\mathbb{R}^{3})\right) \text{ \ with \ }0<r<1,  \label{eq10}
\end{equation}%
then the solution remains smooth on $\left( 0,T\right] $.\ Therefore,
\begin{equation*}
\left( u,\omega \right) \in L^{\infty }\left( 0,T,H^{1}\left( \mathbb{R}%
^{3}\right) \right) \cap L^{2}\left( 0,T,H^{2}\left( \mathbb{R}^{3}\right)
\right) .
\end{equation*}
\end{thm}

For convenience, we will use the following two lemmas will be used in the
proofs of our main results (see, e.g., \cite{Ad, Galdi, Lad}) :

\begin{lem}
Let $\mu ,\lambda $ and $\gamma $ be three parameters that satisfy
\begin{equation*}
1\leq \alpha ,\lambda <\infty ,\text{ \ \ }\frac{1}{\lambda }+\frac{2}{%
\alpha }>1\text{ \ and \ }1+\frac{3}{\gamma }=\frac{1}{\lambda }+\frac{2}{%
\alpha }.
\end{equation*}%
Assume that $f\in H^{1}(\mathbb{R}^{3})$, $\partial _{1}f,\partial _{2}f\in
L^{\alpha }(\mathbb{R}^{3})$ and $\partial _{3}f\in L^{\lambda }(\mathbb{R}%
^{3})$. Then there exists a constant $C=C(\alpha ,\lambda )$ such that
\begin{equation}
\left\Vert f\right\Vert _{L^{\gamma }}\leq C\left\Vert \partial
_{1}f\right\Vert _{L^{\alpha }}^{\frac{1}{3}}\left\Vert \partial
_{2}f\right\Vert _{L^{\alpha }}^{\frac{1}{3}}\left\Vert \partial
_{3}f\right\Vert _{L^{\lambda }}^{\frac{1}{3}},\text{ \ }1\leq \gamma
<\infty .  \label{eq09}
\end{equation}
\end{lem}

\begin{lem}
Let $2\leq \beta \leq 6$ and assume that $f\in H^{1}(\mathbb{R}^{3})$. Then
there exists a constant $C=C(\beta )$ such that%
\begin{eqnarray*}
\left\Vert f\right\Vert _{L^{\beta }} &\leq &C\left\Vert f\right\Vert
_{L^{2}}^{\frac{6-\beta }{2\beta }}\left\Vert \partial _{1}f\right\Vert
_{L^{2}}^{\frac{\beta -2}{2\beta }}\left\Vert \partial _{2}f\right\Vert
_{L^{2}}^{\frac{\beta -2}{2\beta }}\left\Vert \partial _{3}f\right\Vert
_{L^{2}}^{\frac{\beta -2}{2\beta }} \\
&\leq &C\left\Vert f\right\Vert _{L^{2}(\mathbb{R}^{3})}^{\frac{6-\beta }{%
2\beta }}\left\Vert f\right\Vert _{\overset{.}{H}^{1}(\mathbb{R}^{3})}^{%
\frac{3(\beta -2)}{2\beta }}.
\end{eqnarray*}
\end{lem}

Now we are in the position to prove Theorem \ref{th1}.

\begin{pf}
We differentiate the first and the second equation in (\ref{eq1}) with
respect to $x_{3}$, we take the scalar product with $\partial _{3}u$ and $%
\partial _{3}\omega $, respectively and integrate over $\mathbb{R}^{3}$,\ we
get%
\begin{equation}
\frac{1}{2}\frac{d}{dt}\left\Vert \partial _{3}u\right\Vert
_{L^{2}}^{2}+\left\Vert \nabla \partial _{3}u\right\Vert
_{L^{2}}^{2}=-\int\limits_{\mathbb{R}^{3}}\left( \partial _{3}u.\nabla
\right) u.\partial _{3}udx+\int\limits_{\mathbb{R}^{3}}\partial _{3}(\nabla
\times \omega ).\partial _{3}udx.  \label{eq3.1}
\end{equation}%
and%
\begin{eqnarray}
&&\frac{1}{2}\frac{d}{dt}\left\Vert \partial _{3}\omega \right\Vert
_{L^{2}}^{2}+\left\Vert \nabla \partial _{3}\omega \right\Vert
_{L^{2}}^{2}+\left\Vert \nabla .(\partial _{3}\omega )\right\Vert
_{L^{2}}^{2}  \label{eq3.2} \\
&\leq &-\int\limits_{\mathbb{R}^{3}}\left( \partial _{3}u.\nabla \right)
\omega .\partial _{3}\omega dx-2\left\Vert \partial _{3}\omega \right\Vert
_{L^{2}}^{2}+\int\limits_{\mathbb{R}^{3}}\partial _{3}(\nabla \times
u).\partial _{3}\omega dx.  \notag
\end{eqnarray}%
Now, combining (\ref{eq3.1}) and (\ref{eq3.2}), one has after suitable
integration by parts (recall that $\nabla .u=0$)
\begin{eqnarray}
&&\frac{1}{2}\frac{d}{dt}\left[ \left\Vert \partial _{3}u\right\Vert
_{L^{2}}^{2}+\left\Vert \partial _{3}\omega \right\Vert _{L^{2}}^{2}\right]
+\left\Vert \nabla \partial _{3}u\right\Vert _{L^{2}}^{2}+\left\Vert \nabla
\partial _{3}\omega \right\Vert _{L^{2}}^{2}  \notag \\
&\leq &\int\limits_{\mathbb{R}^{3}}\partial _{3}(\nabla \times \omega
).\partial _{3}udx+\int\limits_{\mathbb{R}^{3}}\partial _{3}(\nabla \times
u).\partial _{3}\omega dx-2\left\Vert \partial _{3}\omega \right\Vert
_{L^{2}}^{2}  \label{eq3.3} \\
&&-\int\limits_{\mathbb{R}^{3}}\left( \partial _{3}u.\nabla \right)
u.\partial _{3}udx-\int\limits_{\mathbb{R}^{3}}\left( \partial _{3}u.\nabla
\right) \omega .\partial _{3}\omega dx  \notag \\
&=&A_{1}+A_{2}+A_{3}+A_{4}+A_{5}.  \notag
\end{eqnarray}%
Integrating by parts and using H\"{o}lder's inequality and Young's
inequality (as in \cite{GR_2015_CAMWA}), we derive the estimation of the first three terms on the
right-hand side of (\ref{eq3.3}) as%
\begin{eqnarray*}
A_{1}+A_{2}+A_{3} &=&\int\limits_{\mathbb{R}^{3}}\partial _{3}(\nabla \times
\omega ).\partial _{3}udx+\int\limits_{\mathbb{R}^{3}}\partial _{3}(\nabla
\times u).\partial _{3}\omega dx-2\left\Vert \partial _{3}\omega \right\Vert
_{L^{2}}^{2} \\
&\leq &2\left\Vert \partial _{3}\omega \right\Vert _{L^{2}}^{2}+\frac{1}{2}%
\left\Vert \nabla \partial _{3}u\right\Vert _{L^{2}}^{2}-2\left\Vert
\partial _{3}\omega \right\Vert _{L^{2}}^{2}=\frac{1}{2}\left\Vert \nabla
\partial _{3}u\right\Vert _{L^{2}}^{2}.
\end{eqnarray*}%
For $A_{4}$, using Lemma 2.3 together with the H\"{o}lder inequality and
Young inequality, we find
\begin{eqnarray}
\left\vert A_{4}\right\vert &=&\left\vert \int\limits_{\mathbb{R}^{3}}\left(
\partial _{3}u.\nabla \right) u.\partial _{3}udx\right\vert  \notag \\
&\leq &\left\Vert \partial _{3}u.\partial _{3}u\right\Vert
_{L^{2}}\left\Vert \nabla u\right\Vert _{L^{2}}  \label{eq3.4} \\
&\leq &\left\Vert \partial _{3}u\right\Vert _{\overset{.}{\mathcal{M}}_{2,%
\frac{3}{r}}}\left\Vert \partial _{3}u\right\Vert _{\overset{.}{B}%
_{2,1}^{r}}\left\Vert \nabla u\right\Vert _{L^{2}}  \notag \\
&\leq &\left\Vert \partial _{3}u\right\Vert _{\overset{.}{\mathcal{M}}_{2,%
\frac{3}{r}}}\left\Vert \nabla \partial _{3}u\right\Vert
_{L^{2}}^{r}\left\Vert \partial _{3}u\right\Vert _{L^{2}}^{1-r}\left\Vert
\nabla u\right\Vert _{L^{2}}  \notag
\end{eqnarray}%
by using the following bilinear estimate (see \cite{G1, G2, Lem1}):%
\begin{equation*}
\left\Vert fg\right\Vert _{L^{2}}\leq C\left\Vert f\right\Vert _{\overset{.}{%
\mathcal{M}}_{2,\frac{3}{r}}}\left\Vert g\right\Vert _{\overset{.}{B}%
_{2,1}^{r}}
\end{equation*}%
and the following interpolation inequality \cite{MO} :
\begin{equation*}
\left\Vert w\right\Vert _{\overset{.}{B}_{2,1}^{r}}\leq C\left\Vert
w\right\Vert _{L^{2}}^{1-r}\left\Vert \nabla w\right\Vert _{L^{2}}^{r}.
\end{equation*}%
Similarly, we can bound%
\begin{eqnarray}
\left\vert A_{5}\right\vert &=&\left\vert \int\limits_{\mathbb{R}^{3}}\left(
\partial _{3}u.\nabla \right) \omega .\partial _{3}\omega dx\right\vert
\notag \\
&\leq &\left\Vert \partial _{3}u\right\Vert _{\overset{.}{\mathcal{M}}_{2,%
\frac{3}{r}}}\left\Vert \partial _{3}\omega \right\Vert _{\overset{.}{B}%
_{2,1}^{r}}\left\Vert \nabla \omega \right\Vert _{L^{2}}  \label{eq 3.5} \\
&\leq &\left\Vert \partial _{3}u\right\Vert _{\overset{.}{\mathcal{M}}_{2,%
\frac{3}{r}}}\left\Vert \partial _{3}\omega \right\Vert
_{L^{2}}^{1-r}\left\Vert \nabla \partial _{3}\omega \right\Vert
_{L^{2}}^{r}\left\Vert \nabla \omega \right\Vert _{L^{2}}.  \notag
\end{eqnarray}%
From the above inequalities and (\ref{eq3.3}), we obtain%
\begin{eqnarray*}
&&\frac{1}{2}\frac{d}{dt}\left[ \left\Vert \partial _{3}u\right\Vert
_{L^{2}}^{2}+\left\Vert \partial _{3}\omega \right\Vert _{L^{2}}^{2}\right] +%
\frac{1}{2}\left\Vert \nabla \partial _{3}u\right\Vert
_{L^{2}}^{2}+\left\Vert \nabla \partial _{3}\omega \right\Vert _{L^{2}}^{2}
\\
&\leq &\left\Vert \partial _{3}u\right\Vert _{\overset{.}{\mathcal{M}}_{2,%
\frac{3}{r}}}\left\Vert \nabla \partial _{3}u\right\Vert
_{L^{2}}^{r}\left\Vert \partial _{3}u\right\Vert _{L^{2}}^{1-r}\left\Vert
\nabla u\right\Vert _{L^{2}}+\left\Vert \partial _{3}u\right\Vert _{\overset{%
.}{\mathcal{M}}_{2,\frac{3}{r}}}\left\Vert \partial _{3}\omega \right\Vert
_{L^{2}}^{1-r}\left\Vert \nabla \partial _{3}\omega \right\Vert
_{L^{2}}^{r}\left\Vert \nabla \omega \right\Vert _{L^{2}}
\end{eqnarray*}%
By Young' s inequality $\left( a^{\alpha }b^{1-\alpha }\leq \alpha
a+(1-\alpha )b\leq a+b\text{ with }a,b\geq 0\text{ and }0\leq \alpha \leq
1\right) $, we find%
\begin{eqnarray*}
&&\frac{1}{2}\frac{d}{dt}\left[ \left\Vert \partial _{3}u\right\Vert
_{L^{2}}^{2}+\left\Vert \partial _{3}\omega \right\Vert _{L^{2}}^{2}\right] +%
\frac{1}{2}\left\Vert \nabla \partial _{3}u\right\Vert
_{L^{2}}^{2}+\left\Vert \nabla \partial _{3}\omega \right\Vert _{L^{2}}^{2}
\\
&\leq &\left( \left\Vert \partial _{3}u\right\Vert _{\overset{.}{\mathcal{M}}%
_{2,\frac{3}{r}}}^{\frac{2}{2-r}}\left\Vert \partial _{3}u\right\Vert
_{L^{2}}^{2(\frac{1-r}{2-r})}\left\Vert \nabla u\right\Vert _{L^{2}}^{\frac{2%
}{2-r}}\right) ^{\frac{2-r}{2}}\left( \left\Vert \nabla \partial
_{3}u\right\Vert _{L^{2}}^{2}\right) ^{\frac{r}{2}} \\
&&+3\left( \left\Vert \partial _{3}u\right\Vert _{\overset{.}{\mathcal{M}}%
_{2,\frac{3}{r}}}^{\frac{2}{2-r}}\left\Vert \partial _{3}\omega \right\Vert
_{L^{2}}^{2(\frac{1-r}{2-r})}\left\Vert \nabla \omega \right\Vert _{L^{2}}^{%
\frac{2}{2-r}}\right) ^{\frac{2-r}{2}}\left( \left\Vert \nabla \partial
_{3}\omega \right\Vert _{L^{2}}^{2}\right) ^{\frac{r}{2}} \\
&\leq &C\left\Vert \partial _{3}u\right\Vert _{\overset{.}{\mathcal{M}}_{2,%
\frac{3}{r}}}^{\frac{2}{2-r}}\left\Vert \partial _{3}u\right\Vert
_{L^{2}}^{2(\frac{1-r}{2-r})}\left\Vert \nabla u\right\Vert _{L^{2}}^{\frac{2%
}{2-r}}+C\left\Vert \partial _{3}u\right\Vert _{\overset{.}{\mathcal{M}}_{2,%
\frac{3}{r}}}^{\frac{2}{2-r}}\left\Vert \partial _{3}\omega \right\Vert
_{L^{2}}^{2(\frac{1-r}{2-r})}\left\Vert \nabla \omega \right\Vert _{L^{2}}^{%
\frac{2}{2-r}} \\
&&+\frac{1}{2}\left\Vert \nabla \partial _{3}\omega \right\Vert _{L^{2}}^{2}+%
\frac{1}{2}\left\Vert \nabla \partial _{3}u\right\Vert _{L^{2}}^{2} \\
&=&\frac{1}{2}\left\Vert \nabla \partial _{3}\omega \right\Vert _{L^{2}}^{2}+%
\frac{1}{2}\left\Vert \nabla \partial _{3}u\right\Vert
_{L^{2}}^{2}+C\left\Vert \partial _{3}u\right\Vert _{L^{2}}^{2(\frac{1-r}{2-r%
})}\left( \left( \left\Vert \partial _{3}u\right\Vert _{\overset{.}{\mathcal{%
M}}_{2,\frac{3}{r}}}^{\frac{2}{1-r}}\right) ^{\frac{1-r}{2-r}}\left(
\left\Vert \nabla u\right\Vert _{L^{2}}^{2}\right) ^{\frac{1}{2-r}}\right) \\
&&+C\left\Vert \partial _{3}\omega \right\Vert _{L^{2}}^{2(\frac{1-r}{2-r}%
)}\left( \left( \left\Vert \partial _{3}u\right\Vert _{\overset{.}{\mathcal{M%
}}_{2,\frac{3}{r}}}^{\frac{2}{1-r}}\right) ^{\frac{1-r}{2-r}}\left(
\left\Vert \nabla \omega \right\Vert _{L^{2}}^{2}\right) ^{\frac{1}{2-r}%
}\right) \\
&\leq &\frac{1}{2}\left\Vert \nabla \partial _{3}\omega \right\Vert
_{L^{2}}^{2}+\frac{1}{2}\left\Vert \nabla \partial _{3}u\right\Vert
_{L^{2}}^{2}+C\left\Vert \partial _{3}u\right\Vert _{L^{2}}^{2(\frac{1-r}{2-r%
})}\left( \left\Vert \partial _{3}u\right\Vert _{\overset{.}{\mathcal{M}}_{2,%
\frac{3}{r}}}^{\frac{2}{1-r}}+\left\Vert \nabla u\right\Vert
_{L^{2}}^{2}\right) \\
&&+C\left\Vert \partial _{3}\omega \right\Vert _{L^{2}}^{2(\frac{1-r}{2-r}%
)}\left( \left\Vert \partial _{3}u\right\Vert _{\overset{.}{\mathcal{M}}_{2,%
\frac{3}{r}}}^{\frac{2}{1-r}}+\left\Vert \nabla \omega \right\Vert
_{L^{2}}^{2}\right) ,
\end{eqnarray*}%
which implies that%
\begin{eqnarray*}
&&\frac{1}{2}\frac{d}{dt}(1+\left\Vert \partial _{3}u\right\Vert
_{L^{2}}^{2}+\left\Vert \partial _{3}\omega \right\Vert
_{L^{2}}^{2})+\left\Vert \nabla \partial _{3}u\right\Vert
_{L^{2}}^{2}+\left\Vert \nabla \partial _{3}\omega \right\Vert _{L^{2}}^{2}
\\
&\leq &C(1+\left\Vert \partial _{3}u\right\Vert _{L^{2}}^{2})\left(
\left\Vert \partial _{3}u\right\Vert _{\overset{.}{\mathcal{M}}_{2,\frac{3}{r%
}}}^{\frac{2}{1-r}}+\left\Vert \nabla u\right\Vert _{L^{2}}^{2}\right)
+C(1+\left\Vert \partial _{3}\omega \right\Vert _{L^{2}}^{2})\left(
\left\Vert \partial _{3}u\right\Vert _{\overset{.}{\mathcal{M}}_{2,\frac{3}{r%
}}}^{\frac{2}{1-r}}+\left\Vert \nabla \omega \right\Vert _{L^{2}}^{2}\right)
\\
&\leq &C(1+\left\Vert \partial _{3}u\right\Vert _{L^{2}}^{2}+\left\Vert
\partial _{3}\omega \right\Vert _{L^{2}}^{2})\left( \left\Vert \partial
_{3}u\right\Vert _{\overset{.}{\mathcal{M}}_{2,\frac{3}{r}}}^{\frac{2}{1-r}%
}+\left\Vert \nabla u\right\Vert _{L^{2}}^{2}+\left\Vert \nabla \omega
\right\Vert _{L^{2}}^{2}\right) ,
\end{eqnarray*}%
since $(\frac{1-r}{2-r})<1$. It follows from Gronwall's inequality together
with the energy inequality (\ref{eq1.9}) that%
\begin{eqnarray*}
&&(1+\left\Vert \partial _{3}u(t,.)\right\Vert _{L^{2}}^{2}+\left\Vert
\partial _{3}\omega (t,.)\right\Vert _{L^{2}}^{2}) \\
&\leq &(1+\left\Vert \partial _{3}u_{0}\right\Vert _{L^{2}}^{2}+\left\Vert
\partial _{3}\omega _{0}\right\Vert _{L^{2}}^{2})\exp \left(
C\int_{0}^{t}\left\Vert \partial _{3}u(s,.)\right\Vert _{\overset{.}{%
\mathcal{M}}_{2,\frac{3}{r}}}^{\frac{2}{1-r}}+\left\Vert \nabla
u(s,.)\right\Vert _{L^{2}}^{2}+\left\Vert \nabla \omega (s,.)\right\Vert
_{L^{2}}^{2}ds\right) \\
&\leq &(1+\left\Vert \partial _{3}u_{0}\right\Vert _{L^{2}}^{2}+\left\Vert
\partial _{3}\omega _{0}\right\Vert _{L^{2}}^{2})\exp \left(
C\int_{0}^{t}\left\Vert \partial _{3}u(s,.)\right\Vert _{\overset{.}{%
\mathcal{M}}_{2,\frac{3}{r}}}^{\frac{2}{1-r}}ds+C\left\Vert u_{0}\right\Vert
_{L^{2}}^{2}+\left\Vert \omega _{0}\right\Vert _{L^{2}}^{2}\right) \\
&=&(1+\left\Vert \partial _{3}u_{0}\right\Vert _{L^{2}}^{2}+\left\Vert
\omega _{0}\right\Vert _{L^{2}}^{2})e^{C(\left\Vert u_{0}\right\Vert
_{L^{2}}^{2}+\left\Vert \omega _{0}\right\Vert _{L^{2}}^{2})}\exp \left(
C\int_{0}^{t}\left\Vert \partial _{3}u(s,.)\right\Vert _{\overset{.}{%
\mathcal{M}}_{2,\frac{3}{r}}}^{\frac{2}{1-r}}ds\right)
\end{eqnarray*}%
and \
\begin{equation}
\int_{0}^{t}(\left\Vert \nabla \partial _{3}u(s,.)\right\Vert
_{L^{2}}^{2}+\left\Vert \nabla \partial _{3}\omega (s,.)\right\Vert
_{L^{2}}^{2})ds\leq C.  \label{eq2.14}
\end{equation}%
Here $C$ denotes a constant dependent on the initial data and $\left\Vert
\partial _{3}u(s,.)\right\Vert _{L^{\frac{2}{1-r}}(0,T;\overset{.}{\mathcal{M%
}}_{2,\frac{3}{r}}(\mathbb{R}^{3}))}$.

Now we establish
\begin{equation*}
(u,\omega )\in L^{\infty }(0,T;H^{1})\cap L^{2}(0,T;H^{2}).
\end{equation*}%
Taking the inner product of the equation (\ref{eq1}) with $-\Delta u$ \ and $%
-\Delta \omega $ in $L^{2}(\mathbb{R}^{3})$, respectively, after suitable
integration by parts, by the same calculation as that in \cite{BRTZ}, \cite{G1}, \cite{Jia}, we
obtain for $t\in \left( 0,T\right) $,%
\begin{eqnarray*}
\frac{1}{2}\frac{d}{dt}\left\Vert \nabla u(t)\right\Vert
_{L^{2}}^{2}+\left\Vert \Delta u(t)\right\Vert _{L^{2}}^{2} &=&\int_{\mathbb{%
R}^{3}}(u.\nabla )u.\Delta udx-\int_{\mathbb{R}^{3}}(\nabla \times \omega
).\Delta udx \\
&=&-
\sum_{k=1}^{3}\int_{\mathbb{R}^{3}}\partial _{k}u\cdot (\partial
_{k}u.\nabla u)dx-\int_{\mathbb{R}^{3}}(\nabla \times u).\Delta \omega dx,
\end{eqnarray*}%
\begin{eqnarray*}
&&\frac{1}{2}\frac{d}{dt}\left\Vert \nabla \omega (t)\right\Vert
_{L^{2}}^{2}+\left\Vert \Delta \omega (t)\right\Vert _{L^{2}}^{2}+\left\Vert
\nabla \text{div}\omega (t)\right\Vert _{L^{2}}^{2}+2\left\Vert \nabla
\omega (t)\right\Vert _{L^{2}}^{2} \\
&=&\int_{\mathbb{R}^{3}}(u.\nabla )\omega .\Delta \omega dx-\int_{\mathbb{R}%
^{3}}(\nabla \times u).\Delta \omega dx \\
&=&-
\sum\limits_{k=1}^{3}\int_{\mathbb{R}^{3}}\partial _{k}\omega \cdot
(\partial _{k}u\cdot \nabla \omega )dx-\int_{\mathbb{R}^{3}}(\nabla \times
u).\Delta \omega dx,
\end{eqnarray*}%
where we have used
\begin{equation*}
\int_{\mathbb{R}^{3}}(\nabla \times \omega ).\Delta udx=\int_{\mathbb{R}%
^{3}}(\nabla \times u).\Delta \omega dx.
\end{equation*}%
We sum the above equations to obtain%
\begin{eqnarray*}
&&\frac{1}{2}\frac{d}{dt}(\left\Vert \nabla u(t)\right\Vert
_{L^{2}}^{2}+\left\Vert \nabla \omega (t)\right\Vert
_{L^{2}}^{2})+\left\Vert \Delta u(t)\right\Vert _{L^{2}}^{2}+\left\Vert
\Delta \omega (t)\right\Vert _{L^{2}}^{2} \\
&&+\left\Vert \nabla \text{div}\omega (t)\right\Vert
_{L^{2}}^{2}+2\left\Vert \nabla \omega (t)\right\Vert _{L^{2}}^{2} \\
&\leq &C\left\Vert \nabla u\right\Vert _{L^{3}}^{3}+\left\Vert \nabla
u\right\Vert _{L^{3}}\left\Vert \nabla \omega \right\Vert
_{L^{3}}^{2}+2\left\Vert \nabla u\right\Vert _{L^{2}}\left\Vert \Delta
\omega \right\Vert _{L^{2}} \\
&\leq &C\left\Vert \nabla u\right\Vert _{L^{3}}^{3}+\left( \left\Vert \nabla
u\right\Vert _{L^{3}}^{3}\right) ^{\frac{1}{3}}\left( \left\Vert \nabla
\omega \right\Vert _{L^{3}}^{3}\right) ^{\frac{2}{3}}+C\left\Vert \nabla
u\right\Vert _{L^{2}}^{2}+\frac{1}{4}\left\Vert \Delta \omega \right\Vert
_{L^{2}}^{2} \\
&\leq &C\left\Vert \nabla u\right\Vert _{L^{3}}^{3}+C\left\Vert \nabla
\omega \right\Vert _{L^{3}}^{3}+\frac{1}{4}\left\Vert \Delta \omega
\right\Vert _{L^{2}}^{2}+C\left\Vert \nabla u\right\Vert _{L^{2}}^{2} \\
&\leq &C\left\Vert \nabla u\right\Vert _{L^{2}}^{\frac{3}{2}}\left\Vert
\nabla \partial _{1}u\right\Vert _{L^{2}}^{\frac{1}{2}}\left\Vert \nabla
\partial _{2}u\right\Vert _{L^{2}}^{\frac{1}{2}}\left\Vert \nabla \partial
_{3}u\right\Vert _{L^{2}}^{\frac{1}{2}} \\
&&+C\left\Vert \nabla \omega \right\Vert _{L^{2}}^{\frac{3}{2}}\left\Vert
\nabla \partial _{1}\omega \right\Vert _{L^{2}}^{\frac{1}{2}}\left\Vert
\nabla \partial _{2}\omega \right\Vert _{L^{2}}^{\frac{1}{2}}\left\Vert
\nabla \partial _{3}\omega \right\Vert _{L^{2}}^{\frac{1}{2}}+\frac{1}{4}%
\left\Vert \Delta \omega \right\Vert _{L^{2}}^{2}+C\left\Vert \nabla
u\right\Vert _{L^{2}}^{2} \\
&\leq &C\left\Vert \nabla u\right\Vert _{L^{2}}^{\frac{3}{2}}\left\Vert
\nabla ^{2}u\right\Vert _{L^{2}}\left\Vert \nabla \partial _{3}u\right\Vert
_{L^{2}}^{\frac{1}{2}}+C\left\Vert \nabla \omega \right\Vert _{L^{2}}^{\frac{%
3}{2}}\left\Vert \nabla ^{2}\omega \right\Vert _{L^{2}}\left\Vert \nabla
\partial _{3}\omega \right\Vert _{L^{2}}^{\frac{1}{2}} \\
&&+\frac{1}{4}\left\Vert \Delta \omega \right\Vert _{L^{2}}^{2}+C\left\Vert
\nabla u\right\Vert _{L^{2}}^{2} \\
&=&\left( \left\Vert \nabla ^{2}u\right\Vert _{L^{2}}^{2}\right) ^{\frac{1}{2%
}}\left( C\left\Vert \nabla u\right\Vert _{L^{2}}^{3}\left\Vert \nabla
\partial _{3}u\right\Vert _{L^{2}}\right) ^{\frac{1}{2}}+\left( \left\Vert
\nabla ^{2}\omega \right\Vert _{L^{2}}^{2}\right) ^{\frac{1}{2}}\left(
C\left\Vert \nabla \omega \right\Vert _{L^{2}}^{3}\left\Vert \nabla \partial
_{3}\omega \right\Vert _{L^{2}}\right) ^{\frac{1}{2}} \\
&&+\frac{1}{4}\left\Vert \Delta \omega \right\Vert _{L^{2}}^{2}+C\left\Vert
\nabla u\right\Vert _{L^{2}}^{2} \\
&\leq &\frac{1}{2}\left\Vert \Delta u\right\Vert _{L^{2}}^{2}+C\left\Vert
\nabla u\right\Vert _{L^{2}}^{3}\left\Vert \nabla \partial _{3}u\right\Vert
_{L^{2}}+\frac{1}{4}\left\Vert \Delta \omega \right\Vert
_{L^{2}}^{2}+C\left\Vert \nabla \omega \right\Vert _{L^{2}}^{3}\left\Vert
\nabla \partial _{3}\omega \right\Vert _{L^{2}} \\
&&+\frac{1}{4}\left\Vert \Delta \omega \right\Vert _{L^{2}}^{2}+C\left\Vert
\nabla u\right\Vert _{L^{2}}^{2} \\
&=&\frac{1}{2}(\left\Vert \Delta u\right\Vert _{L^{2}}^{2}+\left\Vert \Delta
\omega \right\Vert _{L^{2}}^{2})+C\left\Vert \nabla u\right\Vert
_{L^{2}}^{2}\left( \left\Vert \nabla u\right\Vert _{L^{2}}\left\Vert \nabla
\partial _{3}u\right\Vert _{L^{2}}\right) \\
&&+C\left\Vert \nabla \omega \right\Vert _{L^{2}}^{2}\left( \left\Vert
\nabla \omega \right\Vert _{L^{2}}\left\Vert \nabla \partial _{3}\omega
\right\Vert _{L^{2}}\right) +\frac{1}{2}(\left\Vert \nabla \omega
\right\Vert _{L^{2}}^{2}+\left\Vert \nabla u\right\Vert _{L^{2}}^{2}) \\
&\leq &\frac{1}{2}(\left\Vert \Delta u\right\Vert _{L^{2}}^{2}+\left\Vert
\Delta \omega \right\Vert _{L^{2}}^{2})+C\left\Vert \nabla u\right\Vert
_{L^{2}}^{2}\left( \left\Vert \nabla u\right\Vert _{L^{2}}^{2}+\left\Vert
\nabla \partial _{3}u\right\Vert _{L^{2}}^{2}\right) \\
&&+C\left\Vert \nabla \omega \right\Vert _{L^{2}}^{2}\left( \left\Vert
\nabla \omega \right\Vert _{L^{2}}^{2}+\left\Vert \nabla \partial _{3}\omega
\right\Vert _{L^{2}}^{2}\right) ,
\end{eqnarray*}%
by using H\"{o}lder inequality and applying (\ref{eq09}) with $\alpha
=\lambda =2$ and $\gamma =6$:
\begin{equation*}
\left\Vert f\right\Vert _{L^{6}}\leq C\left\Vert \partial _{1}f\right\Vert
_{L^{2}}^{\frac{1}{3}}\left\Vert \partial _{2}f\right\Vert _{L^{2}}^{\frac{1%
}{3}}\left\Vert \partial _{3}f\right\Vert _{L^{2}}^{\frac{1}{3}}.
\end{equation*}%
Hence%
\begin{eqnarray*}
&&\frac{d}{dt}(\left\Vert \nabla u(t)\right\Vert _{L^{2}}^{2}+\left\Vert
\nabla \omega (t)\right\Vert _{L^{2}}^{2})+\left\Vert \Delta u(t)\right\Vert
_{L^{2}}^{2}+\left\Vert \Delta \omega (t)\right\Vert _{L^{2}}^{2} \\
&\leq &C(1+\left\Vert \nabla u\right\Vert _{L^{2}}^{2}+\left\Vert \nabla
\partial _{3}u\right\Vert _{L^{2}}^{2}+\left\Vert \nabla \omega \right\Vert
_{L^{2}}^{2}+\left\Vert \nabla \partial _{3}\omega \right\Vert
_{L^{2}}^{2})\left( \left\Vert \nabla u\right\Vert _{L^{2}}^{2}+\left\Vert
\nabla \omega \right\Vert _{L^{2}}^{2}\right)
\end{eqnarray*}%
Using Gronwall's inequality, the energy inequality (\ref{eq1.9}) and the
estimate (\ref{eq2.14}), we conclude that%
\begin{eqnarray*}
&&\left\Vert \nabla u(t,.)\right\Vert _{L^{2}}^{2}+\left\Vert \nabla \omega
(t,.)\right\Vert _{L^{2}}^{2}+\int_{0}^{t}(\left\Vert \Delta
u(s,.))\right\Vert _{L^{2}}^{2}+\left\Vert \Delta \omega (s,.))\right\Vert
_{L^{2}}^{2})ds \\
&\leq &(\left\Vert \nabla u_{0}\right\Vert _{L^{2}}^{2}+\left\Vert \nabla
\omega _{0}\right\Vert _{L^{2}}^{2})\exp \left( C\int_{0}^{t}\left\Vert
\nabla u(s.,)\right\Vert _{L^{2}}^{2}+\left\Vert \nabla \partial
_{3}u(s,.)\right\Vert _{L^{2}}^{2}ds\right) \\
&&\times \exp \left( C\int_{0}^{t}\left\Vert \nabla \omega (s.,)\right\Vert
_{L^{2}}^{2}+\left\Vert \nabla \partial _{3}\omega (s,.)\right\Vert
_{L^{2}}^{2}ds\right) \\
&\leq &C.
\end{eqnarray*}%
for all $0\leq t<T$. Hence%
\begin{equation*}
(u,\omega )\in L^{\infty }(0,T;H^{1}(\mathbb{R}^{3}))\cap L^{2}(0,T;H^{2}(%
\mathbb{R}^{3})),
\end{equation*}%
which gives that $u$ and $\omega $ are smooth. This completes the proof of
Theorem \ref{th1}.
\end{pf}

\begin{re}
Theorem \ref{th1} is still true for the Navier-Stokes equation with $\omega
\equiv 0$, so we give an extension of Serrin's regularity criterion for the
Navier-Stokes equations \cite{Liu2014}.
\end{re}

\section{Acknowledgements}

The authors thank the referees for their invaluable comments and suggestions
which helped improve the paper greatly. This work was done, while the first
author was visiting Catania University in Italy. He thanks Department of
Mathematics and Computer Science at the Catania university for his hospitality.

\end{document}